\documentclass{amsart}

\usepackage[utf8]{inputenc}
\usepackage[T1]{fontenc}
\usepackage{color}

\usepackage{bm}
\usepackage{mathrsfs}
\usepackage{amsmath}
\usepackage{mathtools}
\usepackage{amssymb}
\usepackage{amsthm}
\usepackage{amsfonts}
\usepackage{enumerate}
\usepackage{url}
\usepackage{hyperref}

\usepackage{ulem}

\newcommand{\cb}[1]{\textcolor{black}{#1}}

\theoremstyle{plain}
\newtheorem{thm}{Theorem}[section]
\newtheorem{lem}[thm]{Lemma}

\newtheorem{cor}[thm]{Corollary}

\theoremstyle{definition}
\newtheorem{eg}[thm]{Example}
\newtheorem{rem}[thm]{Remark}
\newtheorem{defn}[thm]{Definition}

\numberwithin{equation}{section}

\newcommand{\ZZ}{{\mathbb Z}}
\newcommand{\QQ}{{\mathbb Q}}
\newcommand{\CC}{{\mathbb C}}

\newcommand{\Li}{\operatorname{Li}}

\begin{document}

\title[FMZV and $p$-MZSV]
{Finite multiple harmonic sums and lower-truncated $p$-adic multiple zeta-star values}

\author{Ku-Yu Fan}
\address{Graduate School of Mathematics, Nagoya University, Furo-cho, Chikusa-ku, Nagoya, 464-8602, Japan.}
\email{ku-yu.fan.d2@math.nagoya-u.ac.jp}

\author{Hidekazu Furusho}
\address{Graduate School of Mathematics, Nagoya University, Furo-cho, Chikusa-ku, Nagoya, 464-8602, Japan.}
\email{furusho@math.nagoya-u.ac.jp}

\date{November 12, 2026.}

\begin{abstract}
We introduce lower-truncated $p$-adic multiple zeta-star values and construct
their harmonic regularizations. We show that they are independent of the branch
of the $p$-adic logarithm and satisfy harmonic product relations. We express
finite multiple harmonic sums as alternating convolutions of harmonic-regularized
$p$-adic multiple zeta values and lower-truncated $p$-adic multiple zeta-star
values. This gives a $p$-adic analogue of the classical truncation identity,
which reduces to the inversion relation between ordinary and star
harmonic-regularized values when the truncation point is $1$.
For truncation point $p$, it relates finite multiple zeta values to $p$-adic
multiple zeta values.
\end{abstract}


\maketitle

\tableofcontents

\section{Introduction}\label{section Introduction}

Multiple zeta values are defined, for positive integers
$k_1,\ldots,k_r$ with $k_r>1$, by
\begin{equation}\label{eq introductory MZV}
  \zeta(k_1,\ldots,k_r)
  =\sum_{0<n_1<\cdots<n_r}
    \frac{1}{n_1^{k_1}\cdots n_r^{k_r}}.
\end{equation}
In depth one, these are special values of the Riemann zeta function.
Their study in higher depths goes back to Euler's work on double
zeta values, and their modern systematic study was developed by
Hoffman~\cite{Hoffman} and Zagier~\cite{Zagier}.
Their series and iterated integral representations give rise to
harmonic and shuffle product relations, respectively.

For an integer $N\geq1$, consider the finite multiple harmonic sum
and the lower-truncated multiple zeta-star function
\[
  \begin{aligned}
    \zeta^{<N}(s_1,\ldots,s_r)
    &:=\sum_{0<n_1<\cdots<n_r<N}
       \frac{1}{n_1^{s_1}\cdots n_r^{s_r}},\\
    \zeta^{\star,N\leq}(s_1,\ldots,s_r)
    &:=\sum_{N\leq n_1\leq\cdots\leq n_r}
       \frac{1}{n_1^{s_1}\cdots n_r^{s_r}}.
  \end{aligned}
\]
The first is entire on $\CC^r$.
The second is a multiple star analogue of the Hurwitz zeta function,
defined by its series in the domain of absolute convergence and
elsewhere by meromorphic continuation.
Writing $\zeta(s_1,\ldots,s_r)$ for the meromorphic continuation
of the multiple zeta series, one has the classical truncation identity
\begin{equation}\label{eq introductory classical identity}
  \zeta^{<N}(s_1,\ldots,s_r)
  =\sum_{j=0}^{r}(-1)^{r-j}
    \zeta(s_1,\ldots,s_j)
    \zeta^{\star,N\leq}(s_r,\ldots,s_{j+1});
\end{equation}
see \cite[Theorem~1.2]{Komori}.
Empty-index factors are understood to be $1$.
The corresponding identity also appears in
\cite[Proposition~2.1]{INY}, in the study of interpolants of
truncated multiple zeta functions.

The purpose of this paper is to construct lower-truncated
$p$-adic multiple zeta-star values and to establish a $p$-adic
analogue of \eqref{eq introductory classical identity}.
For a prime $p$, the second-named  author~\cite{p-adicMZV} defines $p$-adic multiple
zeta values $\zeta_p(\bm{k})$ at admissible indices
$\bm{k}=(k_1,\ldots,k_r)$, with $k_r>1$, by restricted limits
at $z=1$ of multiple polylogarithms analytically continued by
Coleman integration.
These values belong to $\QQ_p$ and are independent of the branch
of the $p$-adic logarithm.
Harmonic regularization extends them to polynomials
$\zeta_{p,*}(\bm{k};T)\in\QQ_p[T]$ for arbitrary positive integer
indices, normalized by $\zeta_{p,*}(1;T)=T$; see~\cite{FJ}, and \cite{IKZ} for the
classical case.
For admissible indices these polynomials are constant and equal
to $\zeta_p(\bm{k})$.

We apply harmonic regularization to analytically continued
lower-truncated multiple star\cb{-}polylogarithms to define
\[
  \zeta_{p,*}^{\star,N\leq}(\bm{k};T)\in\QQ_p[T].
\]
Our first main result, Theorem~\ref{thm lower truncated star properties},
gives an explicit formula for these polynomials in terms of
finite multiple harmonic sums and the untruncated polynomials
$\zeta_{p,*}^{\star}(\bm{k};T)$.
The theorem also establishes independence of the logarithm branch
and the harmonic product relations for star sums.
The ordinary-to-star transformation underlying these relations
is studied in~\cite{IKOO,HoffmanIhara}.
For admissible $\bm{k}$, the polynomials are constant and equal
to restricted limits of the corresponding polylogarithms.
We denote these values by $\zeta_p^{\star,N\leq}(\bm{k})$.
When $N=1$, they recover the usual $p$-adic multiple zeta-star values.

Our second main result, Theorem~\ref{thm Main theorem},
is a $p$-adic analogue of the classical truncation identity
\eqref{eq introductory classical identity}.
For every prime $p$, every integer $N\geq1$, and every index
$(k_1,\ldots,k_r)$ of positive integers, it gives
\begin{equation}\label{eq introductory main regularized}
  \zeta^{<N}(k_1,\ldots,k_r)
  =\sum_{j=0}^{r}(-1)^{r-j}
    \zeta_{p,*}(k_1,\ldots,k_j;T)
    \zeta_{p,*}^{\star,N\leq}(k_r,\ldots,k_{j+1};T)
\end{equation}
in $\QQ_p[T]$.
In particular, the alternating sum on the right is independent
of the regularization parameter $T$.
If $k_1,\ldots,k_r\geq2$, every nonempty index occurring on the
right is admissible, and the identity becomes
\begin{equation}\label{eq introductory main identity}
  \zeta^{<N}(k_1,\ldots,k_r)
  =\sum_{j=0}^{r}(-1)^{r-j}
    \zeta_p(k_1,\ldots,k_j)
    \zeta_p^{\star,N\leq}(k_r,\ldots,k_{j+1})
\end{equation}
in $\QQ_p$; see Corollary~\ref{cor Main theorem admissible factors}.
Thus the same rational finite multiple harmonic sum occurs in
both the classical and the $p$-adic truncation identities.

For $N=1$, the lower truncation disappears and
$\zeta_{p,*}^{\star,1\leq}(\bm{k};T)
=\zeta_{p,*}^{\star}(\bm{k};T)$.
Since $\zeta^{<1}(\bm{k})=0$ for every nonempty index,
\eqref{eq introductory main regularized} reduces to
\begin{equation}\label{eq introductory inversion}
  \sum_{j=0}^{r}(-1)^{r-j}
    \zeta_{p,*}(k_1,\ldots,k_j;T)
    \zeta_{p,*}^{\star}(k_r,\ldots,k_{j+1};T)
  =0
  \qquad(r\geq1).
\end{equation}
This inversion relation follows by applying the harmonic-regularized
evaluation map to the antipode identity for the ordinary harmonic
product, using the explicit antipode formula in
\cite[Theorem~3.2]{H2}.
Our main identity extends this relation to arbitrary lower
truncation points $N$, with the finite multiple harmonic sum
on the left-hand side.

The specialization $N=p$ connects our result with finite multiple
zeta values.
The sums $H_p(\bm{k}):=\zeta^{<p}(\bm{k})$ belong to $\ZZ_p$,
and the family of their residues modulo $p$, as $p$ varies,
defines the finite multiple zeta value associated with $\bm{k}$;
see~\cite{KZ}.
An expression for $H_p(\bm{k})$ in terms of Deligne's $p$-adic
multiple zeta values was conjectured by Akagi, Hirose, and Yasuda
and proved by Jarossay~\cite{Jarossay}; see also
\cite[Theorem~4.1]{Yasuda}.
At $N=p$, the left-hand side of
\eqref{eq introductory main regularized} is $H_p(\bm{k})$.
In Remark~\ref{rem comparison AHY Jarossay}, we compare the
resulting identity with the formula proved by Jarossay.

The proof of our main identity starts from a formal power series
identity for multiple polylogarithms and their truncations.
After Coleman analytic continuation, harmonic regularization
gives \eqref{eq introductory main regularized}.

The paper is organized as follows.
Section~\ref{section Power series} introduces the notation and
proves the two formal truncation identities used in the construction
and in the proof of the main result.
Section~\ref{section lower truncated values} constructs the
lower-truncated $p$-adic multiple zeta-star polynomials and
establishes their branch independence and specialization properties.
Section~\ref{section harmonic products} proves their harmonic
product relations.
Finally, Section~\ref{section Main identity} proves the main
identity and discusses the cases $N=1$ and $N=p$.

\section{Multiple polylogarithms and truncation identities}
\label{section Notation}
\label{section Power series}
We introduce several variants of
the multiple polylogarithms used in this paper and prove two formal truncation identities. 
\subsection{Indices and multiple polylogarithms}

\begin{defn}\label{defn index}\label{defn admissible index}
  An \textit{index} is a finite sequence
  $\bm{k}=(k_1,\ldots,k_r)$ of positive integers, where $r\geq0$.
  The unique index of length zero is denoted by $\emptyset$.
  An index is called \textit{admissible} if it is empty or its last
  entry is greater than $1$.
\end{defn}

Throughout the paper, sums, functions, and regularized polynomials
associated with the empty index are understood to be $1$.
In particular, a factor with an empty sequence of entries or arguments
in any of the formulas below is $1$.

\begin{defn}\label{defn finite MZV}
  For $N\in\ZZ_{\geq1}$ and an index $\bm{k}=(k_1,\ldots,k_r)$,
  the \textit{finite multiple harmonic sum} is
  \[
    \zeta^{<N}(k_1,\ldots,k_r)
    :=
    \sum_{0<n_1<\cdots<n_r<N}
    \frac{1}{n_1^{k_1}\cdots n_r^{k_r}}
    \in\QQ.
  \]
\end{defn}

We use the following four variants of multiple polylogarithms\cb{:}

\begin{defn}\label{defn multivariable MPL}\label{defn star MPL}
\label{defn finite MPL}\label{defn shifted star MPL}
\label{defn one variable MPL}
  Let $\bm{k}=(k_1,\ldots,k_r)$ be a nonempty index and let
  $N\in\ZZ_{\geq1}$. Define
  \begin{align*}
    \Li_{k_1,\ldots,k_r}(z_1,\ldots,z_r)
    &:=
    \sum_{0<n_1<\cdots<n_r}
    \frac{z_1^{n_1}\cdots z_r^{n_r}}
         {n_1^{k_1}\cdots n_r^{k_r}},\\
    \Li^\star_{k_1,\ldots,k_r}(z_1,\ldots,z_r)
    &:=
    \sum_{1\leq n_1\leq\cdots\leq n_r}
    \frac{z_1^{n_1}\cdots z_r^{n_r}}
         {n_1^{k_1}\cdots n_r^{k_r}},\\
    \Li^{<N}_{k_1,\ldots,k_r}(z_1,\ldots,z_r)
    &:=
    \sum_{0<n_1<\cdots<n_r<N}
    \frac{z_1^{n_1}\cdots z_r^{n_r}}
         {n_1^{k_1}\cdots n_r^{k_r}},\\
    \Li^{\star,N\leq}_{k_1,\ldots,k_r}(z_1,\ldots,z_r)
    &:=
    \sum_{N\leq n_1\leq\cdots\leq n_r}
    \frac{z_1^{n_1}\cdots z_r^{n_r}}
         {n_1^{k_1}\cdots n_r^{k_r}}.
  \end{align*}
  These are formal power series over $\QQ$; the third is a polynomial.
  The last two are called the \textit{upper-truncated multiple
  polylogarithm} and the \textit{lower-truncated multiple
  star-polylogarithm}, respectively.
\end{defn}

We write $\bm{z}=(z_1,\ldots,z_r)$ when abbreviating the
arguments. The corresponding one-variable multiple polylogarithm is
\begin{equation}\label{eq one variable specialization}
  \Li_{k_1,\ldots,k_r}(z)
  :=\Li_{k_1,\ldots,k_r}(1,\ldots,1,z)
  =\sum_{0<n_1<\cdots<n_r}
    \frac{z^{n_r}}{n_1^{k_1}\cdots n_r^{k_r}}.
\end{equation}
Evaluation of the finite polynomial gives
\begin{equation}\label{eq finite specialization}
  \Li^{<N}_{\bm{k}}(1,\ldots,1)=\zeta^{<N}(\bm{k}).
\end{equation}

Separating the cases in which consecutive summation variables are
equal gives the following expansion.

\begin{lem}\label{lem star expansion}
  For every nonempty index $\bm{k}=(k_1,\ldots,k_r)$,
  \begin{equation}\label{eq star block expansion}
    \Li^\star_{k_1,\ldots,k_r}(z_1,\ldots,z_r)
    =
    \sum_{\substack{1\leq q\leq r\\
                    0=i_0<i_1<\cdots<i_q=r}}
    \Li_{K_1,\ldots,K_q}(Z_1,\ldots,Z_q),
  \end{equation}
  where
  \[
    K_h=\sum_{\nu=i_{h-1}+1}^{i_h}k_\nu,
    \qquad
    Z_h=\prod_{\nu=i_{h-1}+1}^{i_h}z_\nu
    \qquad(1\leq h\leq q).
  \]
\end{lem}

\begin{proof}
  In a weakly increasing sequence $n_1,\ldots,n_r$, each distinct
  value occurs at consecutive positions. Let $i_h$ be the last
  position occupied by the $h$-th distinct value, so that
  \[
    n_{i_{h-1}+1}=\cdots=n_{i_h}
    \qquad(1\leq h\leq q),
  \]
  with strict inequalities between successive distinct values.
  Equal summation variables add the corresponding denominator
  exponents and multiply the corresponding arguments.
  Summing over all choices of $i_1,\ldots,i_q$ gives the formula.
\end{proof}

\subsection{Two truncation identities}

The first identity expresses a lower-truncated star-polylogarithm
in terms of ordinary star-polylogarithms and finite polynomials.
It will also be used to define its analytic continuation.

\begin{lem}\label{lem shifted star decomposition}
  For every $N\in\ZZ_{\geq1}$ and every index,
  \begin{multline}\label{eq shifted star decomposition}
    \Li^{\star,N\leq}_{k_1,\ldots,k_r}(z_1,\ldots,z_r)
    =\sum_{i=0}^{r}(-1)^i
      \Li^{<N}_{k_i,\ldots,k_1}(z_i,\ldots,z_1)
      \Li^\star_{k_{i+1},\ldots,k_r}(z_{i+1},\ldots,z_r)
  \end{multline}
  in $\QQ[[z_1,\ldots,z_r]]$.
\end{lem}

\begin{proof}
  The empty-index case follows from our conventions.
  For $r\geq1$, put
  \[
    \cb{A}_i:=\Li^{<N}_{k_i,\ldots,k_1}(z_i,\ldots,z_1)
         \Li^\star_{k_{i+1},\ldots,k_r}(z_{i+1},\ldots,z_r)
    \qquad(0\leq i\leq r)
  \]
  and
  \[
    R_i:=
    \sum_{\substack{n_1,\ldots,n_r\geq1\\
                    N>n_1>\cdots>n_{i+1}\\
                    n_{i+1}\leq\cdots\leq n_r}}
    \frac{z_1^{n_1}\cdots z_r^{n_r}}
         {n_1^{k_1}\cdots n_r^{k_r}}
    \qquad(0\leq i<r).
  \]
  For $1\leq i<r$, splitting the domain of $\cb{A}_i$ according to
  $n_i\leq n_{i+1}$ or $n_i>n_{i+1}$ gives
  \[
    \cb{A}_i=R_{i-1}+R_i.
  \]
  Moreover, $R_{r-1}=\cb{A}_r$ and
  $\Li^{\star,N\leq}_{\bm{k}}(\bm{z})=\cb{A}_0-R_0$.
  Successive substitution therefore gives
  \[
    \cb{A}_0-R_0=\sum_{i=0}^{r}(-1)^i\cb{A}_i.
  \]
  All decompositions are coefficientwise identities of formal
  power series, proving the assertion.
\end{proof}

The second identity is a polylogarithmic version of the
classical truncation formula in
\cite[Theorem~1.2]{Komori}.
The corresponding classical identity also appears in
\cite[Proposition~2.1]{INY}, in the study of interpolants
of truncated multiple harmonic sums.
We will use its polylogarithmic version to prove
the main theorem.

\begin{lem}\label{lem power series identity}
  For every $N\in\ZZ_{\geq1}$ and every index,
  \begin{multline}\label{eq MPL identity}
    \Li^{<N}_{k_1,\ldots,k_r}(z_1,\ldots,z_r)
    =\sum_{j=0}^{r}(-1)^{r-j}
      \Li_{k_1,\ldots,k_j}(z_1,\ldots,z_j)
      \Li^{\star,N\leq}_{k_r,\ldots,k_{j+1}}
        (z_r,\ldots,z_{j+1})
  \end{multline}
  in $\QQ[[z_1,\ldots,z_r]]$.
\end{lem}

\begin{proof}
  Again, the empty-index case is immediate. For $r\geq1$, put
  \[
    \cb{B}_j:=\Li_{k_1,\ldots,k_j}(z_1,\ldots,z_j)
         \Li^{\star,N\leq}_{k_r,\ldots,k_{j+1}}
           (z_r,\ldots,z_{j+1})
    \qquad(0\leq j\leq r)
  \]
  and
  \[
    S_j:=
    \sum_{\substack{n_1,\ldots,n_r\geq1\\
                    n_1<\cdots<n_{j+1}\\
                    n_{j+1}\geq\cdots\geq n_r\geq N}}
    \frac{z_1^{n_1}\cdots z_r^{n_r}}
         {n_1^{k_1}\cdots n_r^{k_r}}
    \qquad(0\leq j<r).
  \]
  For $1\leq j<r$, the alternatives $n_j<n_{j+1}$ and
  $n_j\geq n_{j+1}$ give $\cb{B}_j=S_j+S_{j-1}$.
  Since $S_0=\cb{B}_0$ and
  $\Li^{<N}_{\bm{k}}(\bm{z})=\cb{B}_r-S_{r-1}$, we obtain
  \[
    \cb{B}_r-S_{r-1}=\sum_{j=0}^{r}(-1)^{r-j}\cb{B}_j.
  \]
  This is the required formal power series identity.
\end{proof}

\section[Lower-truncated p-adic multiple zeta-star values]
{Lower-truncated $p$-adic multiple zeta-star values}
\label{section lower truncated values}
We construct lower-truncated $p$-adic multiple zeta-star values and
their harmonic regularizations via Coleman analytic continuation.
Theorem~\ref{thm lower truncated star properties} gives an explicit
evaluation formula and establishes branch independence and harmonic
product relations.

\subsection{Analytic continuation and harmonic regularization}

Fix a prime $p$. We work over $\QQ_p$ and the completion $\CC_p$
of its algebraic closure, with $|p|_p=p^{-1}$.
Fix a branch $\log^\varpi$ of the $p$-adic logarithm such that
$\log^\varpi(p)=\varpi\in\QQ_p$.

Coleman integration \cite{C} gives analytic continuations of the
one-variable multiple polylogarithms, denoted by
$\Li^\varpi_{\bm{k}}(z)$; see
\cite[Definition~2.9]{p-adicMZV}.
We also use the multivariable Coleman functions of \cite{BF,FJ},
denoted by $\Li^\varpi_{\bm{k}}(z_1,\ldots,z_r)$.
All identities between these functions are considered on their
common domains. The identity principle ensures that identities
valid in the region of convergence remain valid for the
corresponding Coleman functions; see \cite[\S2]{BF}.

We write $\lim'$ for the restricted $p$-adic limit of
\cite[Notation~2.12]{p-adicMZV}. Thus the limit is required to
have the same value along every sequence
$\{t_n\}\subset\CC_p\setminus\{1\}$ approaching $1$
for which $\QQ_p(t_1,t_2,\ldots)/\QQ_p$ has finite absolute
ramification index.

\begin{defn}\label{defn p-adic MZV}
  For a nonempty admissible index $\bm{k}$, the
  \textit{$p$-adic multiple zeta value} is
  \[
    \zeta_p(\bm{k})
    :=\lim_{z\to1}{}'\Li^\varpi_{\bm{k}}(z).
  \]
\end{defn}

This limit exists, belongs to $\QQ_p$, and is independent of
$\varpi$; see \cite[Theorems~2.13, 2.18 and 2.25]{p-adicMZV}.

For arbitrary indices, we use harmonic regularization.
For a function with a local logarithmic expansion
\[
  F(t)=\sum_{m=0}^{M}f_m^\varpi(t)
       \bigl(-\log^\varpi(1-t)\bigr)^m,
\]
where each coefficient is analytic over $\QQ_p$ at $t=1$, set
\begin{equation}\label{eq regularization operation}
  \operatorname{Reg}_{t=1}^{\varpi}F
  :=\sum_{m=0}^{M}f_m^\varpi(1)T^m.
\end{equation}
The local expansions and corresponding tangential morphisms described in \cite[Proposition~4.1 and \S6]{FJ} apply to the functions used below, including finite sums and products and functions obtained by substituting products of arguments as in \eqref{eq star block expansion}.
Our convention for the formal variable $T$ has the opposite sign
to that in \cite{FJ}.
The operation $\operatorname{Reg}_{t=1}^{\varpi}$ is
multiplicative, since evaluation of the analytic coefficients at
$t=1$ is multiplicative; thus, for functions $F$ and $G$ admitting such expansions, we have
\begin{equation}\label{eq regularization multiplicativity}
  \operatorname{Reg}_{t=1}^{\varpi}(FG)
  =\bigl(\operatorname{Reg}_{t=1}^{\varpi}F\bigr)
   \bigl(\operatorname{Reg}_{t=1}^{\varpi}G\bigr).
\end{equation}
In particular,
\begin{equation}\label{eq regularization polynomial}
  \operatorname{Reg}_{t=1}^{\varpi}P(t,\ldots,t)
  =P(1,\ldots,1)
  \qquad\bigl(P\in\QQ_p[z_1,\ldots,z_r]\bigr).
\end{equation}

\begin{defn}\label{defn ordinary harmonic regularization}
  For any index $\bm{k}=(k_1,\ldots,k_r)$, write
  \[
    \zeta_{p,*}(\bm{k};T)
    :=\operatorname{Reg}_{t=1}^{\varpi}
      \Li^\varpi_{\bm{k}}(t,\ldots,t),
    \qquad
    \zeta_{p,*}(\bm{k}):=\zeta_{p,*}(\bm{k};0).
  \]
\end{defn}

These are the usual harmonic-regularized $p$-adic multiple zeta
polynomials, the $p$-adic counterparts of the harmonic-regularized multiple zeta values
$\zeta_*(\boldsymbol{k};T)$ introduced in \cite{IKZ}.
They satisfy $\zeta_{p,*}(1;T)=T$ and are independent
of $\varpi$ by \cite[Theorem~6.4]{FJ}.
For admissible $\bm{k}$, they are constant and equal to
$\zeta_p(\bm{k})$ by \cite[Corollary~5.3]{FJ}.

Star values are expressed as sums of ordinary values
over the indices obtained by replacing any selection
of commas by addition.
The algebraic properties of this transformation,
including its compatibility with harmonic products,
are studied in \cite{IKOO,HoffmanIhara}.
We use this relation to define the harmonic-regularized
$p$-adic multiple zeta-star polynomials.

\begin{defn}\label{defn p-adic star MZV}
  For a nonempty index $\bm{k}=(k_1,\ldots,k_r)$, define
  \[
    \zeta_{p,*}^{\star}(\bm{k};T)
    :=\sum_{\substack{1\leq q\leq r\\
                       0=i_0<i_1<\cdots<i_q=r}}
       \zeta_{p,*}(K_1,\ldots,K_q;T),
  \]
  where $K_h=\sum_{\nu=i_{h-1}+1}^{i_h}k_\nu$.
  Set $\zeta_{p,*}^{\star}(\bm{k})
  :=\zeta_{p,*}^{\star}(\bm{k};0)$.
  If $\bm{k}$ is admissible, this polynomial is constant;
  its value is denoted by $\zeta_p^\star(\bm{k})$ and called
  the \textit{$p$-adic multiple zeta-star value}.
\end{defn}

Indeed, admissibility gives $K_q\geq k_r>1$ in every summand.
Branch independence follows from that of the ordinary
regularized polynomials.

When we restrict \eqref{eq star block expansion} to
$z_1=\cdots=z_r=t$, its arguments become positive powers of $t$.
The following specialization rule explains how these powers
affect regularization.

\begin{lem}\label{lem weighted strict specialization}
  Let $\bm{\ell}=(\ell_1,\ldots,\ell_s)$ be a nonempty index,
  and let $d_1,\ldots,d_s$ be positive integers satisfying
  \begin{equation}\label{eq trailing ones condition}
    d_i=1
    \qquad\text{whenever }
    \ell_i=\ell_{i+1}=\cdots=\ell_s=1.
  \end{equation}
  Then
  \begin{align}
    \operatorname{Reg}_{t=1}^{\varpi}
    \Li^\varpi_{\bm{\ell}}(t^{d_1},\ldots,t^{d_s})
    &=\zeta_{p,*}(\bm{\ell};T),
    \label{eq weighted strict specialization}\\
    \operatorname{Reg}_{t=1}^{\varpi}
    \Li^{\star,\varpi}_{\bm{\ell}}(t^{d_1},\ldots,t^{d_s})
    &=\zeta_{p,*}^{\star}(\bm{\ell};T).
    \label{eq weighted star specialization}
  \end{align}
\end{lem}

Here $\Li^{\star,\varpi}$ is defined by replacing each ordinary
multiple polylogarithm in \eqref{eq star block expansion}
by its Coleman analytic continuation $\Li^{\varpi}$.
For $p$-adic (twisted) multiple star-polylogarithms constructed by
Coleman integration, see \cite{FKMT}.

\begin{proof}
  Put $e_0=0$ and $e_i=d_1+\cdots+d_i$ for $1\leq i\leq s$,
  and set
  \[
    y_i:=\prod_{\nu=e_{i-1}+1}^{e_i}x_\nu.
  \]
  Setting all $x_\nu=t$ gives $y_i=t^{d_i}$.
If $\ell_s>1$, apply \cite[Corollary~5.3]{FJ} to the
products $y_i$ and then set all $x_\nu=t$. This gives
\[
  \operatorname{Reg}_{t=1}^{\varpi}
  \Li^\varpi_{\bm{\ell}}(t^{d_1},\ldots,t^{d_s})
  =\zeta_p(\bm{\ell})
  =\zeta_{p,*}(\bm{\ell};T),
\]
proving \eqref{eq weighted strict specialization}.  If all $\ell_i=1$, condition
  \eqref{eq trailing ones condition} forces all $d_i=1$, and
  the assertion is the definition of
  $\zeta_{p,*}(\bm{\ell};T)$.

  In the remaining case, let $j<s$ be the position of the last
  entry greater than $1$. Then $d_{j+1}=\cdots=d_s=1$.
  Proposition~5.1 of \cite{FJ}, applied to the arguments $y_i$
  and followed by the diagonal restriction, gives
  \begin{align*}
    &\operatorname{Reg}_{t=1}^{\varpi}
      \Li^\varpi_{\bm{\ell}}(t^{d_1},\ldots,t^{d_s})
    =
      \operatorname{Reg}_{t=1}^{\varpi}
      \Li^\varpi_{\bm{\ell}}
      (\underbrace{1,\ldots,1}_{j},
       \underbrace{t,\ldots,t}_{s-j}).
  \end{align*}
  The same proposition with all $d_i=1$ identifies the
  right-hand side with $\zeta_{p,*}(\bm{\ell};T)$.
  This is the comparison used in the proof of
  \cite[Lemma~6.2]{FJ}, and proves
  \eqref{eq weighted strict specialization}.

  For the star version, apply \eqref{eq star block expansion}.
  The entries and arguments in each summand are
  \[
    K_h=\sum_{\nu=i_{h-1}+1}^{i_h}\ell_\nu,
    \qquad
    t^{D_h},
    \qquad
    D_h=\sum_{\nu=i_{h-1}+1}^{i_h}d_\nu.
  \]
  If $K_h=\cdots=K_q=1$, positivity of the original entries
  implies that each of these sums contains exactly one entry
  from the terminal string of $1$'s of $\bm{\ell}$.
  Hence $D_h=\cdots=D_q=1$.
  Formula \eqref{eq weighted strict specialization} therefore
  applies to every summand, and their sum is
  $\zeta_{p,*}^\star(\bm{\ell};T)$ by definition.
\end{proof}

\subsection{Lower-truncated star values and harmonic products}
\label{section harmonic products}

Define $\Li^{\star,\varpi,N\leq}_{\bm{k}}$ by replacing
the multiple star-polylogarithms on the right-hand side of
\eqref{eq shifted star decomposition} by their analytically
continued functions. The finite factors are polynomials, so
this gives an analytic continuation of
$\Li^{\star,N\leq}_{\bm{k}}$ as a Coleman function.
The expansion results of \cite[Proposition~4.3 and \S6]{FJ}
give a local logarithmic expansion of its diagonal restriction
of the form used in \eqref{eq regularization operation}.

\begin{defn}\label{defn regularized lower truncated pMZSV}
  For $N\in\ZZ_{\geq1}$ and any index $\bm{k}$, define the
  \textit{harmonic-regularized lower-truncated $p$-adic multiple
  zeta-star polynomial} by
  \[
    \zeta_{p,*}^{\star,N\leq}(\bm{k};T)
    :=\operatorname{Reg}_{t=1}^{\varpi}
       \Li^{\star,\varpi,N\leq}_{\bm{k}}(t,\ldots,t)
    \in\QQ_p[T],
  \]
  and set
  $\zeta_{p,*}^{\star,N\leq}(\bm{k})
   :=\zeta_{p,*}^{\star,N\leq}(\bm{k};0)$.
\end{defn}

For example,
\begin{equation}\label{eq lower truncated normalization}
  \zeta_{p,*}^{\star,N\leq}(1;T)
  =T-\sum_{n=1}^{N-1}\frac{1}{n}.
\end{equation}

We first recall the product on indices that will be used in the
theorem below.

Let $\mathcal I$ be the $\QQ$-vector space with basis consisting
of all indices, including $\emptyset$.
For $a\in\ZZ_{>0}$ and an index $\bm{u}=(u_1,\ldots,u_m)$,
write
\[
  (a,\bm{u})=(a,u_1,\ldots,u_m),
  \qquad (a,\emptyset)=(a).
\]
For a finite linear combination
$X=\sum_{\bm{w}}c_{\bm{w}}\bm{w}\in\mathcal I$, set
\[
  (a,X):=\sum_{\bm{w}}c_{\bm{w}}(a,\bm{w}).
\]

\begin{defn}\label{defn star harmonic product}
  The \textit{star harmonic product} is the $\QQ$-bilinear map
  \[
    *_\star:\mathcal I\times\mathcal I\longrightarrow\mathcal I
  \]
  determined recursively by
  \[
    \emptyset *_\star\bm{u}
    =\bm{u} *_\star\emptyset=\bm{u}
  \]
  and
  \begin{align*}
    (a,\bm{u}) *_\star(b,\bm{v})
    =\bigl(a,\bm{u} *_\star(b,\bm{v})\bigr)
       +\bigl(b,(a,\bm{u}) *_\star\bm{v}\bigr)
        -\bigl(a+b,\bm{u} *_\star\bm{v}\bigr),
  \end{align*}
  where $a,b\in\ZZ_{>0}$ and $\bm{u},\bm{v}$ are arbitrary,
  possibly empty, indices.
\end{defn}

This is the harmonic product corresponding to
weakly increasing sums; see \cite{IKOO,HoffmanIhara}.
We extend the assignment
$\bm{k}\mapsto\zeta_{p,*}^{\star,N\leq}(\bm{k};T)$
by $\QQ$-linearity to a map from $\mathcal I$ to $\QQ_p[T]$.

The following theorem gives an explicit evaluation of the
lower-truncated star polynomials and establishes their harmonic
product relations.

\begin{thm}
\label{thm lower truncated star properties}
\label{prop lower truncated star regularization}
\label{prop lower truncated star harmonic}
  Let $N\in\ZZ_{\geq1}$. The following assertions hold.
  \begin{enumerate}
    \item[\textup{(i)}]
      For every index
      $\bm{k}=(k_1,\ldots,k_r)$,
      \begin{equation}\label{eq lower truncated regularization formula}
        \zeta_{p,*}^{\star,N\leq}(\bm{k};T)
        =\sum_{i=0}^{r}(-1)^i
          \zeta^{<N}(k_i,\ldots,k_1)
          \zeta_{p,*}^{\star}(k_{i+1},\ldots,k_r;T).
      \end{equation}
      In particular, this polynomial is independent of $\varpi$.
      If $\bm{k}$ is admissible, it is constant, and
      \[
        \lim_{t\to1}{}'
          \Li^{\star,\varpi,N\leq}_{\bm{k}}(t,\ldots,t)
        =\zeta_{p,*}^{\star,N\leq}(\bm{k})\in\QQ_p.
      \]

    \item[\textup{(ii)}]
      For any indices
      $\bm{u},\bm{v}$,
      \[
        \zeta_{p,*}^{\star,N\leq}(\bm{u};T)
        \zeta_{p,*}^{\star,N\leq}(\bm{v};T)
        =\zeta_{p,*}^{\star,N\leq}(\bm{u} *_\star\bm{v};T).
      \]
      In particular, for $a,b\in\ZZ_{>0}$,
      \begin{align*}
        &\zeta_{p,*}^{\star,N\leq}(a;T)
         \zeta_{p,*}^{\star,N\leq}(b;T)
         \\&\qquad=
         \zeta_{p,*}^{\star,N\leq}(a,b;T)
         +\zeta_{p,*}^{\star,N\leq}(b,a;T)
         -\zeta_{p,*}^{\star,N\leq}(a+b;T).
      \end{align*}
  \end{enumerate}
\end{thm}

\begin{proof}
  \textup{(i)} Let $\bm{\ell}=(\ell_1,\ldots,\ell_s)$ be a nonempty index,
  and let $d_1,\ldots,d_s$ be positive integers satisfying
  \eqref{eq trailing ones condition}.
  Every nonempty terminal sequence
  $(\ell_{i+1},\ldots,\ell_s)$ inherits this condition.
  Apply \eqref{eq shifted star decomposition} with
  $z_i=t^{d_i}$ and then apply
  $\operatorname{Reg}_{t=1}^{\varpi}$.
  Formula \eqref{eq weighted star specialization} evaluates
  every infinite factor, and
  \eqref{eq regularization polynomial} evaluates every finite
  factor at $t=1$. By \eqref{eq regularization multiplicativity}, we obtain
  \begin{align*}
    &\operatorname{Reg}_{t=1}^{\varpi}
      \Li^{\star,\varpi,N\leq}_{\bm{\ell}}
        (t^{d_1},\ldots,t^{d_s})
    =
      \sum_{i=0}^{s}(-1)^i
      \zeta^{<N}(\ell_i,\ldots,\ell_1)
      \zeta_{p,*}^{\star}(\ell_{i+1},\ldots,\ell_s;T).
  \end{align*}
  Taking all $d_i=1$ proves
  \eqref{eq lower truncated regularization formula}.
  For general $d_i$, comparison with this formula gives
  \begin{equation}\label{eq weighted lower star specialization}
    \operatorname{Reg}_{t=1}^{\varpi}
    \Li^{\star,\varpi,N\leq}_{\bm{\ell}}
      (t^{d_1},\ldots,t^{d_s})
    =\zeta_{p,*}^{\star,N\leq}(\bm{\ell};T)
  \end{equation}
  under condition \eqref{eq trailing ones condition}.
  Since every term on the right-hand side of
  \eqref{eq lower truncated regularization formula} is
  independent of $\varpi$, so is its left-hand side.

  If $\bm{k}$ is nonempty and admissible, every nonempty
  terminal sequence in that formula is admissible.
  Hence the polynomial is a constant $c\in\QQ_p$.
  In the local expansion
  \[
    \Li^{\star,\varpi,N\leq}_{\bm{k}}(t,\ldots,t)
    =\sum_{m=0}^{M}f_m^\varpi(t)
      \bigl(-\log^\varpi(1-t)\bigr)^m,
  \]
  we therefore have $f_0^\varpi(1)=c$ and
  $f_m^\varpi(1)=0$ for $m\geq1$.
  By analyticity, there are functions $g_m^\varpi$, analytic at $t=1$, such that
  \[
    \Li^{\star,\varpi,N\leq}_{\bm{k}}(t,\ldots,t)-c
      =(t-1)\sum_{m=0}^{M}g_m^\varpi(t)
        \bigl(-\log^\varpi(1-t)\bigr)^m.
  \]
  Since
  \[
    \lim_{t\to1}{}'(t-1)
      \bigl(\log^\varpi(1-t)\bigr)^m=0
    \qquad(m\geq0)
  \]
  by \cite[Lemma~2.14]{p-adicMZV} (applied to
  $\varepsilon=1-t$), the asserted limit follows.
  The empty-index case follows from our conventions.

  \smallskip
  \noindent\textup{(ii)} The cases involving the empty index are
  immediate. Consider two lower-truncated multiple
  star-polylogarithms in independent sets of variables. In depth one,
  inclusion-exclusion gives
  \[
    \Li^{\star,N\leq}_{a}(x)\Li^{\star,N\leq}_{b}(y)
    =
     \Li^{\star,N\leq}_{a,b}(x,y)
     +\Li^{\star,N\leq}_{b,a}(y,x)
     -\Li^{\star,N\leq}_{a+b}(xy).
  \]
  The two weak orderings cover all pairs of summation variables,
  and their common diagonal must be subtracted once.
  On this diagonal, the denominator exponents add and the
  arguments multiply. Repeating this decomposition gives
  the recursion defining $*_\star$ in arbitrary depths.

  By the identity principle, the product identity remains
  valid for the corresponding Coleman functions.
  Set all original variables equal to $t$.
  In every resulting term, an argument is $t$ if it comes
  from one original argument, or $t^2$ if two arguments have
  been multiplied. The latter accompanies an entry of the
  form $a+b\geq2$. Hence every entry equal to $1$ retains
  the argument $t$, so
  \eqref{eq trailing ones condition} holds for every term.
  Apply $\operatorname{Reg}_{t=1}^{\varpi}$
  using \eqref{eq regularization multiplicativity} and
  \eqref{eq weighted lower star specialization}.
  This proves the product formula.
\end{proof}

\begin{defn}\label{defn lower truncated pMZSV}
\label{defn shifted p-adic star MZV}
  For $N\in\ZZ_{\geq1}$ and an admissible index $\bm{k}$,
  the \textit{lower-truncated $p$-adic multiple zeta-star value}
  is
  \[
    \zeta_p^{\star,N\leq}(\bm{k})
    :=\lim_{t\to1}{}'
       \Li^{\star,\varpi,N\leq}_{\bm{k}}(t,\ldots,t).
  \]
\end{defn}

Theorem~\ref{thm lower truncated star properties}\textup{(i)}
shows that this value belongs to $\QQ_p$, is independent of
$\varpi$, and satisfies
\[
  \zeta_{p,*}^{\star,N\leq}(\bm{k};T)
  =\zeta_p^{\star,N\leq}(\bm{k})
  \qquad\text{for admissible }\bm{k}.
\]
For $N=1$, it agrees with $\zeta_p^\star(\bm{k})$.

\begin{eg}\label{eg lower truncated star product}
  For positive integers $a,b,c$, the star harmonic product gives
  \[
    \begin{aligned}
      &\Li^{\star,N\leq}_{a,b}(x_1,x_2)\,
       \Li^{\star,N\leq}_{c}(y)
       =\Li^{\star,N\leq}_{a,b,c}(x_1,x_2,y)
       +\Li^{\star,N\leq}_{a,c,b}(x_1,y,x_2)\\
      &\quad+\Li^{\star,N\leq}_{c,a,b}(y,x_1,x_2)
       -\Li^{\star,N\leq}_{a,b+c}(x_1,x_2y)
       -\Li^{\star,N\leq}_{a+c,b}(x_1y,x_2).
    \end{aligned}
  \]
  After Coleman analytic continuation, set $x_1=x_2=y=t$.
  The two merged terms have the specializations
  \[
    \begin{aligned}
      &\operatorname{Reg}_{t=1}^{\varpi}
        \Li^{\star,\varpi,N\leq}_{a,b+c}(t,t^2)
        =\zeta_{p,*}^{\star,N\leq}(a,b+c;T),\\
      &\operatorname{Reg}_{t=1}^{\varpi}
        \Li^{\star,\varpi,N\leq}_{a+c,b}(t^2,t)
        =\zeta_{p,*}^{\star,N\leq}(a+c,b;T),
    \end{aligned}
  \]
  by \eqref{eq weighted lower star specialization}:
  the squared arguments accompany the entries $b+c\geq2$ and $a+c\geq2$.
  Using \eqref{eq regularization multiplicativity}, we therefore obtain
  \[
    \begin{aligned}
      &\zeta_{p,*}^{\star,N\leq}(a,b;T)\,
       \zeta_{p,*}^{\star,N\leq}(c;T)
       =\zeta_{p,*}^{\star,N\leq}(a,b,c;T)
       +\zeta_{p,*}^{\star,N\leq}(a,c,b;T)\\
      &\quad+\zeta_{p,*}^{\star,N\leq}(c,a,b;T)
       -\zeta_{p,*}^{\star,N\leq}(a,b+c;T)
       -\zeta_{p,*}^{\star,N\leq}(a+c,b;T).
    \end{aligned}
  \]
\end{eg}

\section{The $p$-adic truncation identity}
\label{section Main identity}

We now apply the harmonic regularization of
Section~\ref{section lower truncated values} to the formal
identity of Lemma~\ref{lem power series identity}.
This gives an identity valid for all positive integer indices.
We then discuss the specializations $N=1$ and $N=p$,
and compare the latter with a formula proved by Jarossay \cite{Jarossay}.

\begin{thm}\label{thm Main theorem}
  Let $N\in\ZZ_{\geq1}$.
  For every $r\in\ZZ_{\geq0}$ and integers
  $k_1,\ldots,k_r\geq1$, the following identity holds in
  $\QQ_p[T]$:
  \begin{equation}\label{eq Main theorem polynomial}
    \zeta^{<N}(k_1,\ldots,k_r)
    =\sum_{j=0}^{r}(-1)^{r-j}
      \zeta_{p,*}(k_1,\ldots,k_j;T)
      \zeta_{p,*}^{\star,N\leq}
        (k_r,\ldots,k_{j+1};T).
  \end{equation}
  In particular, the right-hand side is independent of $T$.
  Evaluating at $T=0$ gives
  \begin{equation}\label{eq Main theorem}
    \zeta^{<N}(k_1,\ldots,k_r)
    =\sum_{j=0}^{r}(-1)^{r-j}
      \zeta_{p,*}(k_1,\ldots,k_j)
      \zeta_{p,*}^{\star,N\leq}(k_r,\ldots,k_{j+1}).
  \end{equation}
\end{thm}

\begin{proof}
  The case $r=0$ follows from the empty-index conventions.
  Suppose that $r\geq1$.
  Analytically continue the infinite factors in
  \eqref{eq MPL identity} as Coleman functions, with a fixed
  logarithm branch $\log^\varpi$.
  By the identity principle recalled in
  Section~\ref{section lower truncated values},
  the resulting functions satisfy the same identity
  on their common domain.

  Set $z_1=\cdots=z_r=t$ and apply
  $\operatorname{Reg}_{t=1}^{\varpi}$.
  The left-hand side is a polynomial in $t$;
  by \eqref{eq regularization polynomial} and
  \eqref{eq finite specialization}, its regularization is
  $\zeta^{<N}(k_1,\ldots,k_r)$.
  By definition, the two infinite factors in each
  summand regularize to the corresponding polynomials
  on the right-hand side of
  \eqref{eq Main theorem polynomial}.
  Multiplicativity of regularization \cb{in \eqref{eq regularization multiplicativity}} proves
  \eqref{eq Main theorem polynomial}.
  Its left-hand side is independent of $T$, and
  evaluation at $T=0$ gives \eqref{eq Main theorem}.
\end{proof}

\begin{cor}\label{cor Main theorem admissible factors}
  Under the hypotheses of Theorem~\ref{thm Main theorem},
  if $k_1,\ldots,k_r\geq2$, then
  \begin{equation}\label{eq Main theorem admissible factors}
    \zeta^{<N}(k_1,\ldots,k_r)
    =\sum_{j=0}^{r}(-1)^{r-j}
      \zeta_p(k_1,\ldots,k_j)
      \zeta_p^{\star,N\leq}(k_r,\ldots,k_{j+1}).
  \end{equation}
\end{cor}

\begin{proof}
  Every nonempty index occurring on the right-hand side
  is admissible.
  Hence each regularized polynomial in
  \eqref{eq Main theorem polynomial} is constant and
  equal to the corresponding unregularized value.
\end{proof}

\begin{rem}
\label{rem Main theorem N one}
  When $N=1$, the lower truncation disappears, so
  \[
    \zeta_{p,*}^{\star,1\leq}(\bm{k};T)
    =\zeta_{p,*}^{\star}(\bm{k};T).
  \]
  For every nonempty index,
  $\zeta^{<1}(\bm{k})=0$.
  Thus \eqref{eq Main theorem polynomial} becomes
  \begin{equation}\label{eq Main theorem N one}
    \sum_{j=0}^{r}(-1)^{r-j}
      \zeta_{p,*}(k_1,\ldots,k_j;T)
      \zeta_{p,*}^{\star}(k_r,\ldots,k_{j+1};T)
    =0
    \qquad(r\geq1).
  \end{equation}
  This is the standard inversion relation between
  ordinary and star harmonic-regularized values;
  see \cite[Theorem~3.2]{H2} for its
  interpretation in terms of the antipode.
  Theorem~\ref{thm Main theorem} extends this relation
  to lower-truncated star values, with the finite
  multiple harmonic sum on the left-hand side.
\end{rem}

\begin{rem}\label{rem:equiv}
Given the inversion relation \eqref{eq Main theorem N one}, 
the equation \eqref{eq Main theorem polynomial} is equivalent to
formula \eqref{eq lower truncated regularization formula}. 
Indeed, substituting \eqref{eq lower truncated regularization formula} into the right-hand side of 
\eqref{eq Main theorem polynomial}
and reordering, the coefficient of $\zeta^{<N}(k_{r-i+1},\dots,k_r)$ becomes
 \[
  \sum_{j=0}^{r-i}(-1)^{r-i-j}\zeta_{p,*}(k_1,\dots,k_j;T)
  \zeta^{\star}_{p,*}(k_{r-i},\dots,k_{j+1};T),
\]
 which vanishes for $i<r$ by  \eqref{eq Main theorem N one}
and equals $1$ for $i=r$. The converse is proved in the same way.
\end{rem}

\begin{rem}
\label{rem finite multiple zeta values}
  The specialization $N=p$ connects the main identity
  with finite multiple zeta values.
  Put
  \[
    H_p(\bm{k}):=\zeta^{<p}(\bm{k}).
  \]
  Since all denominators in its defining sum are
  prime to $p$, we have $H_p(\bm{k})\in\ZZ_p$.
  Recall that the finite multiple zeta value
  associated with $\bm{k}$ is
  \[
    \zeta_{\mathcal A}(\bm{k})
    :=\bigl(H_p(\bm{k})\bmod p\bigr)_p
    \in\mathcal A,
    \qquad
    \mathcal A:=
    \left(\prod_{p\ \mathrm{prime}}\mathbb F_p\right)
    \Big/
    \left(\bigoplus_{p\ \mathrm{prime}}\mathbb F_p\right),
  \]
  where $\mathbb F_p=\ZZ/p\ZZ$;
  see \cite{KZ}.

 For each prime $p$, set $N=p$ in \eqref{eq Main theorem}.
The resulting alternating sum equals $H_p(\bm{k})$
and hence lies in $\ZZ_p$.
Reducing these sums modulo their respective primes and
taking their class in $\mathcal A$ gives
$\zeta_{\mathcal A}(\bm{k})$.
\end{rem}

\begin{rem}
\label{rem comparison AHY Jarossay}
  Another expression for $H_p(\bm{k})$ was conjectured
  by Akagi--Hirose--Yasuda and proved by Jarossay:
  \begin{equation}\label{eq comparison AHY Jarossay}
    H_p(k_1,\ldots,k_r)
    =\sum_{j=0}^{r}(-1)^{k_{j+1}+\cdots+k_r}
      \zeta_p^{\mathrm{De}}(k_1,\ldots,k_j)
      \zeta_p^{\mathrm{De},\downarrow}(k_r,\ldots,k_{j+1}).
  \end{equation}
  See \cite[Equation (0.\cb{4}.8)]{Jarossay} and
  \cite[Theorem~4.1]{Yasuda}.
  Here $\zeta_p^{\mathrm{De}}$ denotes Deligne's
  $p$-adic multiple zeta value (\cite{F07}).
  Our convention for $H_p$ does not include the
  weight factor $p^{k_1+\cdots+k_r}$.
  The down-arrow values are given by the
  $p$-adically convergent series
  \begin{align*}
    &\zeta_p^{\mathrm{De},\downarrow}(k_1,\ldots,k_r)
    =
      \sum_{\ell_1,\ldots,\ell_r\geq0}
      p^{\ell_1+\cdots+\ell_r}
      \prod_{i=1}^{r}\binom{k_i+\ell_i-1}{\ell_i}
      \zeta_p^{\mathrm{De}}(k_1+\ell_1,\ldots,k_r+\ell_r).
  \end{align*}
  The empty-index values are again $1$.

  Both \eqref{eq Main theorem} at $N=p$ and
  \eqref{eq comparison AHY Jarossay} involve an
  initial segment and a reversed terminal segment
  of the index.
  Their signs and their $p$-adic factors differ:
  our formula \eqref{eq Main theorem} \cb{uses} the depth $r-j$ in the sign and
  lower-truncated star values, whereas
  \eqref{eq comparison AHY Jarossay} uses the weight
  $k_{j+1}+\cdots+k_r$ and the down-arrow
  transformation of Deligne's values.
  For example, in depth one their comparison gives
\[
H_p(k)=\zeta_{p\cb{,*}}(k)-\zeta_{p\cb{,*}}^{\star,p\leq}(k),
\qquad
H_p(k)=\zeta_p^{\mathrm{De}}(k)
       +(-1)^k\zeta_p^{\mathrm{De},\downarrow}(k)\cb{.}
\]
  Theorem~\ref{thm Main theorem} follows from harmonic
  regularization of a formal truncation identity and
  holds for every integer $N\geq1$.
  The equation \eqref{eq comparison AHY Jarossay}
  expresses $H_p$ through
  Frobenius values and convergent series of such values.
\end{rem}

\medskip
{\bf Declaration of AI Use.}
During the preparation of this manuscript, ChatGPT (OpenAI), a generative AI tool, was used for English-language editing throughout the manuscript, including improvements to the structure, clarity, and wording of the prose.


\end{document}